\documentclass{amsart}

\usepackage{amsmath}

\usepackage{euscript}
\usepackage{mathrsfs}
\usepackage{bbold}
\usepackage{amssymb}
\usepackage{mathabx}

\usepackage{tikz}
\usetikzlibrary{arrows}

\usepackage{float}

\RequirePackage{color}
\definecolor{myred}{rgb}{0.75,0,0}
\definecolor{mygreen}{rgb}{0,0.5,0}
\definecolor{myblue}{rgb}{0,0,0.65}

\RequirePackage{ifpdf}
\ifpdf
  \IfFileExists{pdfsync.sty}{\RequirePackage{pdfsync}}{}
  \RequirePackage[pdftex,
   colorlinks = true,
   urlcolor = myblue, 
   citecolor = mygreen, 
   linkcolor = myred, 
   pagebackref,
   bookmarksopen=true]{hyperref}
\else
  \RequirePackage[hypertex]{hyperref}
  \fi

  \usepackage{url}

\RequirePackage{ae, aecompl, aeguill} 

\newcommand{\nc}{\newcommand} \newcommand{\renc}{\renewcommand}

    \def\RM{{\mathbb{R}}}

    \def\LC{{\mathcal{L}}}

    \def\RC{{\mathcal{R}}}

\def\XSS{{\mathscr{X}}}

\nc{\todo}[1]{ {\color{red}XXX #1 XXX}}

\def\to{\rightarrow}

\def\longto{\longrightarrow}

\nc{\triright}{\stackrel{[1]}{\to}}
\nc{\longtriright}{\stackrel{[1]}{\longto}}

\nc{\Hb}{H^\bullet}

\nc{\Br}{\mathcal{B}}
\nc{\HotRR}{{}_R\mathcal{K}_R}
\nc{\HotR}{\mathcal{K}_R}
\nc{\excise}[1]{}
\nc{\defect}{\text{df}}
\nc{\h}[1]{\underline{H}_{#1}}

\nc{\Ga}{\mathbb{G}_a} 
\nc{\Gm}{\mathbb{G}_m} 

\nc{\Perv}{{\mathbf{P}}}

\nc{\IH}{{\mathrm{IH}}}

\nc{\ic}{\mathbf{IC}}

\nc{\gl}{{\mathfrak{gl}}}
\renc{\sl}{{\mathfrak{sl}}}
\renc{\sp}{{\mathfrak{sp}}}

\renc{\Im}{\textrm{Im}}

\nc{\HBM}{H^{BM}}

\nc{\St}{\mathrm{St}}
\nc{\rot}{\mathrm{rot}}
\nc{\ext}{\mathrm{ext}}
\nc{\Tilt}{\mathrm{Tilt}}
\nc{\gen}{\mathrm{gen}}
\nc{\Graph}{\mathrm{Graph}}

\nc{\simto}{\stackrel{\sim}{\to}}
\nc{\simfrom}{\stackrel{\sim}{\leftarrow}}

\nc{\gbmod}{\mathrm{-gmod-}}

\nc{\gmod}{\mathrm{-gmod}}

\nc{\Parity}{\mathrm{Parity}}

\nc{\mult}{\mathrm{mult}}

\nc{\Hecke}{\textrm{H}}

\nc{\geom}{\mathrm{geom}}
\nc{\Soe}{\mathrm{Soe}}
\nc{\Abe}{\mathrm{Abe}}
\nc{\diag}{\mathrm{diag}}

\nc{\fin}{\textrm{finite}}

\nc{\reflect}{\RC}
\nc{\Chi}{\XSS}
\nc{\pt}{\mathrm{pt}}
\nc{\odd}{\textrm{odd}}
\nc{\even}{\textrm{even}}

\nc{\weights}{\textrm{weights}}

\nc{\adhoc}{\textrm{adhoc}}
\nc{\fund}{\textrm{fund}}

\newtheorem{thm}{Theorem}[section]

\theoremstyle{definition}

\theoremstyle{remark}
\newtheorem{remark}[thm]{Remark}

\tikzset{%
  every neuron/.style={
    circle,
    minimum size=.1cm
  }
}

\title[]{Is deep learning a useful tool for the pure mathematician?}

\author[]{Geordie Williamson}
\address{University of Sydney,Australia.}
\email{g.williamson@sydney.edu.au}

\begin{document}

\begin{abstract}
  A personal and informal account of what a pure mathematician might
  expect when using tools from deep learning in their research.
\end{abstract}

\maketitle



\section{Introduction} \label{sec:intro}

Over the last decade, deep learning has found
countless applications throughout industry and science. However,
its impact on pure mathematics has been modest. This is perhaps
surprising, as some of the tasks at which deep learning excels---like
playing the board-game Go or finding patterns in complicated
structures---appear to present similar difficulties to problems encountered
in research mathematics. On the other hand, the ability to \emph{reason}---probably the
single most important defining characteristic of mathematical
enquiry---remains a central unsolved problem in artificial
intelligence.  Thus, mathematics can be seen as an important litmus
test as to what modern artificial
intelligence can and cannot do.

There is great potential for interaction between mathematics and
machine learning.\footnote{In 1948, Turing  \cite[\S 6]{turing}
  identifies games, mathematics, cryptography and language translation
  and acquisition as five ``suitable branches of
  thought'' in which experimentation with machine intelligence might be fruitful.}
However, there is also a lot of hype, and it is easy for the
mathematician to be put off. In my experience, it remains hard to use deep
learning to aid my mathematical research. However it is possible. One
also has the sense that the potential, once the right tools have been uncovered, is significant.

This is a very informal survey of what a working
mathematician might expect when using the tools of deep
learning on mathematics problems. I outline some of the
beautiful ideas behind deep learning. I
also give some practical hints for using these tools. 
I finish with some examples where deep learning has been used
productively in pure mathematics research. (I hope it goes
  without saying that the impact of deep learning on applied
  mathematics has been enormous.)

Finally, in my experience, the more one uses the tools of deep
learning, the more difficult it becomes not to ask oneself foundational
questions about why they work. This raises an entirely different set
of questions. Although fascinating, the mathematical theory of deep
learning is not the focus here.

\begin{remark} The elephant in the room in any discussion today of
  deep learning is the recent success of ChatGPT and other large
  language models. The internet is full of examples of ChatGPT
  doing both very well and very poorly on reasoning and mathematics
  problems. It seems likely that large
  language models will be able to interact well with proof assistants
  in the near future (see e.g. \cite{han2021proof,jiang2023draft}). It is also likely that a greater role will be
  played in mathematics research by very large models, possibly with emergent capabilities
  (``foundation models'' in the language of the excellent \cite{foundation}).
  The impacts of such developments on mathematics
  are difficult to predict. In this article I will ignore these
  questions entirely. Thus I will restrict myself to situations in
  which deep learning can be used by mathematicians without access to
  these large models.
\end{remark}

{\small 
\subsection{About the author} I am a pure mathematician, working
mostly in geometric representation theory and related fields. I began an ongoing
collaboration with DeepMind in 2020, on possible interactions of machine
learning and mathematics, and have been fascinated by the subject ever
since.\footnote{My understanding of this landscape has benefitted enormously
from discussions with Charles Blundell, Lars Buesing, Alex Davies,  Joel
Gibson, Georg Gottwald, Camilo Libedinsky, S\'ebastien Racaniere, Carlos Simpson, 
Grzegorz Swirszcz, Petar Veličkovi\'c, Adam Wagner, Théophane Weber
and Greg Yang. Without their help this journey would have been 
slower and much more painful. This article is based on a lecture given
at the Fields Institute Symposium on the future of mathematical
research, which was held at the instigation of Akshay Venkatesh. I
would like to thank the organizers for organising a wonderful and
thought-provoking meeting,
and in particular Maia Fraser for very useful feedback on an earlier
version of this article.} }

\section{What is a neural network?} \label{sec:nn}

Artificial neural networks emulate the biological
neural networks present in the brains of humans and other
animals. Typically this emulation takes place on a computer. The
idea of doing so is very natural. See \cite{mccullochpitts, turing} for
remarkable early accounts.

A cartoon picture of a neuron imagines it as a unit with several
inputs and a single output, which may then be connected to other neurons:
\[
  \begin{tikzpicture}[yscale=.4]
    \node (z) at (0,0) {$\bullet$};
    \node (a) at (-2,2) {};
    \node (b) at (-2,1) {};
    \node (c) at (-2,0) {};
    \node (d) at (-2,-1) {};
    \node (e) at (-2,-2) {};
    \draw[->] (a) to[out=0,in=110] (z);
    \draw[->] (b) to[out=0,in=150] (z);
    \draw[->] (c) to[out=0,in=180] (z); 
    \draw[->] (d) to[out=0,in=210] (z); 
    \draw[->] (e) to[out=0,in=250] (z);
    \node (r1) at (3,1) {};
    \node (r2) at (3,0) {};
    \node (r3) at (3,-1) {};
    \node (rr) at (1.7,0) {};
    \node (r) at (1.4,0) {};
    \draw[-] (z) to (rr);
    \draw[->] (r) to[out=0,in=180] (r1);
    \draw[->] (r) to[out=0,in=180] (r2);
    \draw[->] (r) to[out=0,in=180] (r3);
    \end{tikzpicture}
  \]
Neurons ``fire'' by emitting electrical charge along their axon. We
may encode the charges arriving along each node by a real number,  in
which case the charge emitted by a neuron is given by
\[
  \begin{array}{cc}
  \begin{tikzpicture}[yscale=.6,xscale=1.3]
    \node (z) at (0,0) {$\bullet$};
    \node (r) at (2,0) {};
    \node (a) at (-2,2) {};
    \node (b) at (-2,1) {};
    \node (c) at (-2,0) {};
    \node (d) at (-2,-1) {};
    \node (e) at (-2,-2) {};
    \draw[->] (a) to[out=0,in=110] node[fill,white] {$x_1$} node {$x_1$}  (z);
    \draw[->] (b) to[out=0,in=150] node[fill,white] {$x_2$} node {$x_2$}  (z);
    \draw[->] (c) to[out=0,in=180] node[fill,white] {$x_3$} node {$x_3$}  (z);
    \draw[->] (d) to[out=0,in=210] node[fill,white] {$x_4$} node {$x_4$}  (z);
    \draw[->] (e) to[out=0,in=250] node[fill,white] {$x_5$} node {$x_5$}  (z);
        \node (r1) at (3,1) {};
    \node (r2) at (3,0) {};
    \node (r3) at (3,-1) {};
    \node (rr) at (1.5,0) {$z$};
    \node (r) at (1.5,0) {};
    \draw[-] (z) to (rr);
    \draw[->] (r) to[out=0,in=180] (r1);
    \draw[->] (r) to[out=0,in=180] (r2);
    \draw[->] (r) to[out=0,in=180] (r3);
  \end{tikzpicture} &
    \end{array} \quad   z  = f \left (\sum x_i \right )
  \]
  where $f$ is a (typically monotone increasing and non-linear) \emph{activation
    function}. Soon we will assume that our activation function is
  fixed\footnote{and equal to ``ReLU'': $f(x) = \max(0,x)$}, however at this
  level of precision the reader is encouraged to imagine something
  like $f(x) = \tanh(x)$.  The
  activation function is meant to model the non-linear response curves
  of neurons to stimuli. For example, some neurons may not fire until a
  certain charge is reached at their source.\footnote{In biological neural nets
  there is typically large variation in the responses of neurons to
  stimuli, depending on where they are in the brain (see
  e.g. \cite{receptive}). This is one of the many features of
  biological neural nets 
  that is usually ignored when building artificial neural networks.}

  Another   important feature of neurons is that their firing may be
  excitatory or inhibitory of downstream neurons to varying degrees. In order to
  account for this, one allows modification of the input charges via
  weights (the $w_i$):
  \begin{equation} \label{eq:neuron_weights}
  \begin{array}{cc}
  \begin{tikzpicture}[yscale=.6,xscale=1.3]
    \node (z) at (0,0) {$\bullet$};
    \node (r) at (2,0) {};
    \node (a) at (-2,2) {};
    \node (b) at (-2,1) {};
    \node (c) at (-2,0) {};
    \node (d) at (-2,-1) {};
    \node (e) at (-2,-2) {};
    \draw[->] (a) to[out=0,in=110] node[fill,white] {$x_1$} node {$x_1$}  (z);
    \draw[->] (b) to[out=0,in=150] node[fill,white] {$x_2$} node {$x_2$}  (z);
    \draw[->] (c) to[out=0,in=180] node[fill,white] {$x_3$} node {$x_3$}  (z);
    \draw[->] (d) to[out=0,in=210] node[fill,white] {$x_4$} node {$x_4$}  (z);
    \draw[->] (e) to[out=0,in=250] node[fill,white] {$x_5$} node {$x_5$}  (z);
        \node (r1) at (3,1) {};
    \node (r2) at (3,0) {};
    \node (r3) at (3,-1) {};
    \node (rr) at (1.5,0) {$z$};
    \node (r) at (1.5,0) {};
    \draw[-] (z) to (rr);
    \draw[->] (r) to[out=0,in=180] (r1);
    \draw[->] (r) to[out=0,in=180] (r2);
    \draw[->] (r) to[out=0,in=180] (r3);
  \end{tikzpicture} &
    \end{array} \quad   z  = f \left (\sum w_i x_i \right )
  \end{equation}
Thus positive and negative weights correspond to excitatory and
inhibitory connections respectively.

Having settled on a crude mathematical model of a single neuron, we
may then assemble them together to form a \emph{neural network}:
\[
\begin{tikzpicture}[xscale=.8, yscale=.6,x=1.5cm, y=1.5cm, >=stealth]

\foreach \m/\l [count=\y] in {1,2,3,4,5}
  \node [every neuron/.try, neuron \m/.try] (input-\m) at (0,2.5-\y) {$\bullet$};

\foreach \m [count=\y] in {1,2,3,4}
  \node [every neuron/.try, neuron \m/.try ] (hidden-\m) at (2,2-\y) {$\bullet$};

\foreach \m [count=\y] in {1,2}
  \node [every neuron/.try, neuron \m/.try ] (output-\m) at (4,1-\y) {$\bullet$};

\foreach \l [count=\i] in {1,2,3,4,5}
  \draw [<-] (input-\i) -- ++(-1,0);

\foreach \l [count=\i] in {1,2}
  \draw [->] (output-\i) -- ++(1,0);

\foreach \i in {1,...,4,5}
  \foreach \j in {1,2,3,4}
    \draw [->] (input-\i) -- (hidden-\j);

\foreach \i in {1,2,3,4}
  \foreach \j in {1,2}
    \draw [->] (hidden-\i) -- (output-\j);

\end{tikzpicture}
\]
Implicit in this picture is the
assignment of a weight to each edge. Thus our neural network yields a
function which takes real valued inputs (5 in the above picture), and
outputs real values (2 above), via repeated application of 
\eqref{eq:neuron_weights} at each node.
  
This is a good picture for the layperson to have in mind. It is 
useful to visualize the complex interconnectedness present in
artificial neural networks, as well as the locality of the computation
taking place. However for the mathematician, one can
explain things a little differently. The configuration
\[
\begin{tikzpicture}[xscale=.7, yscale=.5,x=1.5cm, y=1.5cm, >=stealth]

\foreach \m/\l [count=\y] in {1,2,3,4,5}
  \node [every neuron/.try, neuron \m/.try] (input-\m) at (0,2.5-\y) {$\bullet$};

\foreach \m [count=\y] in {1,2,3,4}
  \node [every neuron/.try, neuron \m/.try ] (hidden-\m) at (2,2-\y) {$\bullet$};




\foreach \i in {1,...,4,5}
  \foreach \j in {1,2,3,4}
    \draw [->] (input-\i) -- (hidden-\j);


\end{tikzpicture}
\]
is simply a complicated way of drawing a $5 \times 4$ matrix. In other
words, we can rewrite our neural network above economically in the form
\[
  \RM^5 \stackrel{W_1}{\longto} \RM^4 \stackrel{f}{\longto} \RM^4
  \stackrel{W_2}{\longto} \RM^2 \stackrel{f}{\longto} \RM^2 
\]
where the $W_i$ are linear maps determined by matrices of weights, and $f$ is shorthand for the
coordinatewise application of our activation function $f$.

For the purposes of this article, a \emph{vanilla neural
  network}\footnote{One often encounters the term ``Multi
    Layer Perceptron (MLP)'' in the literature.} is a gadget of the
form
\[
  \RM^{d_1} \stackrel{A_1}{\longto} \RM^{d_2}
  \stackrel{f}{\longto} \RM^{d_2} \stackrel{A_2}{\longto}
  \RM^{d_3}
\stackrel{f}{\longto} \RM^{d_3} \stackrel{A_3}{\longto}  \dots \stackrel{f}{\longto}  
\RM^{d_{\ell-1}} \stackrel{A_{\ell-1}}{\longto}\RM^{d_{\ell}} 
\]
where $A_i$ are affine linear maps. We refer to
$\RM^{d_1}, \RM^{d_2}, \dots, \RM^{d_\ell}$ as the \emph{layers} of the network. In order to
simplify the discussion, we always assume that our activation function
$f$ is given by
ReLU (the ``rectified linear unit''), that is
\[
f\left ( \sum \lambda_i e_i \right ) = \sum \max(\lambda_i, 0) e_i
\]
where the $e_i$ are standard basis vectors.

  \begin{remark} We make the following remarks:
    \begin{enumerate}
    \item The attentive reader might have observed a sleight of hand
      above, where we suddenly allowed affine linear maps in our
      definition of a vanilla neural net. This can be justified as
      follows: In biological neural nets both the charge triggering a
      neuron to fire, as well as the charge emitted, varies across
      the neural network. This suggests that each activation
      function should have parameters, i.e. be given by $x
      \mapsto f(x + a) + b$ for varying $a, b \in \RM$ at each
      node. Things just got a lot more complicated! Affine linear maps
      circumvent this issue: by
      adding the possibility of affine linear maps one gets the same
      degree of expressivity, with a much simpler setup.
    \item We only consider ReLU activation functions below. This one
      of the standard choices, and provides a useful
      simplification. However one shouldn't forget that it is possible
      to vary activation functions. 
\item We have tried to motivate the above discussion of neural
  networks as some imitation of neural activity. It is important
  to keep in mind that this is a very loose metaphor at best. However
  I do find it useful in understanding and motivating basic concepts. For an
  excellent account along these lines by an excellent mathematician,
  the reader is referred to \cite{astonishing}.
\item The alert reader will notice that we have implicitly assumed
  above that our graphs representing neural networks do not have any
  cycles or loops. This is again a simplification, and it is desirable
  in certain situations (e.g. in recurrent neural networks) to
  allow loops.
    \end{enumerate}
  \end{remark}

  Vanilla neural networks are often referred to as \emph{fully-connected}
  because each neuron is connected to every neuron in
  the next layer. This is almost opposite to the situation encountered in the brain,
  where remarkably sparse neural networks are found. The connection
  pattern of neurons is referred to the \emph{architecture} of the
  neural network. As well as vanilla neural networks, important
  artificial neural network architectures 
  include \emph{convolutional neural networks}, \emph{graph neural
    networks} and \emph{transformers}. Constraints of length prohibit
  us from discussing these architectures in any depth.

\begin{remark}
  More generally, nowadays the term ``neural network'' is often used
  to refer to any program in which the output depends in a smooth way
  on the input (and thus the program can be updated via some form of
  gradient descent). We ignore this extra generality here.
\end{remark}

\section{Motivation for deep learning}

In order to understand deep learning, it is useful to keep in mind the
tasks at which it first excelled. One of the most important
such examples is image classification. For example, we might want to
classify hand-written digits:
\[
  \begin{array}{c}\includegraphics[width=.8cm]{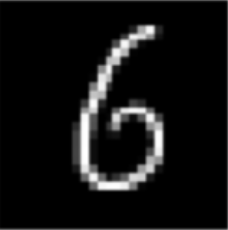}\end{array} \mapsto 6 \quad
  \begin{array}{c}\includegraphics[width=.8cm]{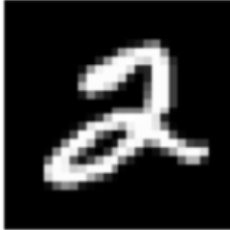}\end{array} \mapsto 2
\]
Here each digit is given as (say) a $28 \times 28$ matrix of
grayscale values between 0 and 255. This is a task which is
effortless for us, but is traditionally difficult for computers.

We can imagine that our brain contains a function which sees a
hand-written digit and produces a probability distribution on $\{ 0, 1,
  \dots, 9 \}$, i.e. ``what digit we think it is''.\footnote{I can convince
  myself that my brain produces a probability distribution and not a
  yes/no answer by recalling my efforts to decipher my grandmother's
  letters when I was a child.} We might attempt to imitate this function with a neural
network.

Let us consider a
simpler problem in which we try to decide whether a
hand-written digit is a 6 or not:
\[
  \begin{array}{c}\includegraphics[width=.8cm]{img/6.png}\end{array} \mapsto \text{``yes"} \quad
  \begin{array}{c}\includegraphics[width=.8cm]{img/2.png}\end{array} \mapsto \text{``no"}
\]
We assume that we have ``training data''
consisting of images labelled by ``6'' or ``not 6''. As a first attempt we might consider a network having a single linear layer:
\[
\RM^{28 \times 28} \stackrel{A}{\longto} \RM \quad 
\stackrel{1/(1+e^{-x})}{\xrightarrow{\hspace*{1.5cm}}} \quad \RM.
\]
Here $A$ is affine linear, and the second function (the ``logistic
function''\footnote{a.k.a. \emph{sigmoid} in the machine learning literature}) is a
convenient way of converting an arbitrary real number into a
probability. Thus, positive values of $A$ mean that we think our image
is a 6, and negative values of $A$ mean we think it is not.

We will be successful if we can find a hyperplane separating all vectors
corresponding to 6's (red dots) from those that do not represent 6's
(blue dots):
\[ \begin{tikzpicture}[xscale=.5,yscale=.5]
    \draw (0,-1.5) to (2,3.5);
    \node[red] (x) at (1,-.7) {$\bullet$};
    \node[red] (x) at (5,.7) {$\bullet$};
    \node[red] (x) at (4,.7) {$\bullet$};
    \node[red] (x) at (4,-1) {$\bullet$};
    \node[blue] (x) at (.2,1.2) {$\bullet$};
    \node[blue] (x) at (-2.5,.5) {$\bullet$};
    \node[blue] (x) at (-3,1.5) {$\bullet$};
    \node[blue] (x) at (1,2) {$\bullet$};
    \node (x) at (10,1) {\tiny $\begin{array}{c} \text{vectors in } \\
                                  \RM^{28 \times 28} \end{array}$};
                              \node (blah) at (2.5,3) {$+$};
                              \node (blah) at (1.3,3.4) {$-$};
                    \end{tikzpicture} \]       
Of course, it may not be possible to find such a hyperplane. Also,
even if we find a hyperplane separating red and blue dots, it is not
clear that such a rule would \emph{generalize}, to correctly predict
whether an unseen image (i.e. image not in our training data)
represents a 6 or not. Remarkably, techniques of this form (for
example \emph{logistic regression}, \emph{Support Vector Machines
  (SVMs)}, \dots) \emph{do} work in many simple
learning scenarios. Given training data (e.g. a large set of
vectors labelled with ``yes'' and ``no'') the optimal separating
hyperplane may be found easily.\footnote{For a striking mathematical example 
  of support vector machines see \cite{he2022learning}, where SVMs are
  trained to 
  distinguish simple and non-simple finite groups, by inspection of
  their multiplication table.}

\section{What is deep learning?}

In many classification problems the classes are not linearly
separable:
\[
  \includegraphics[width=5cm]{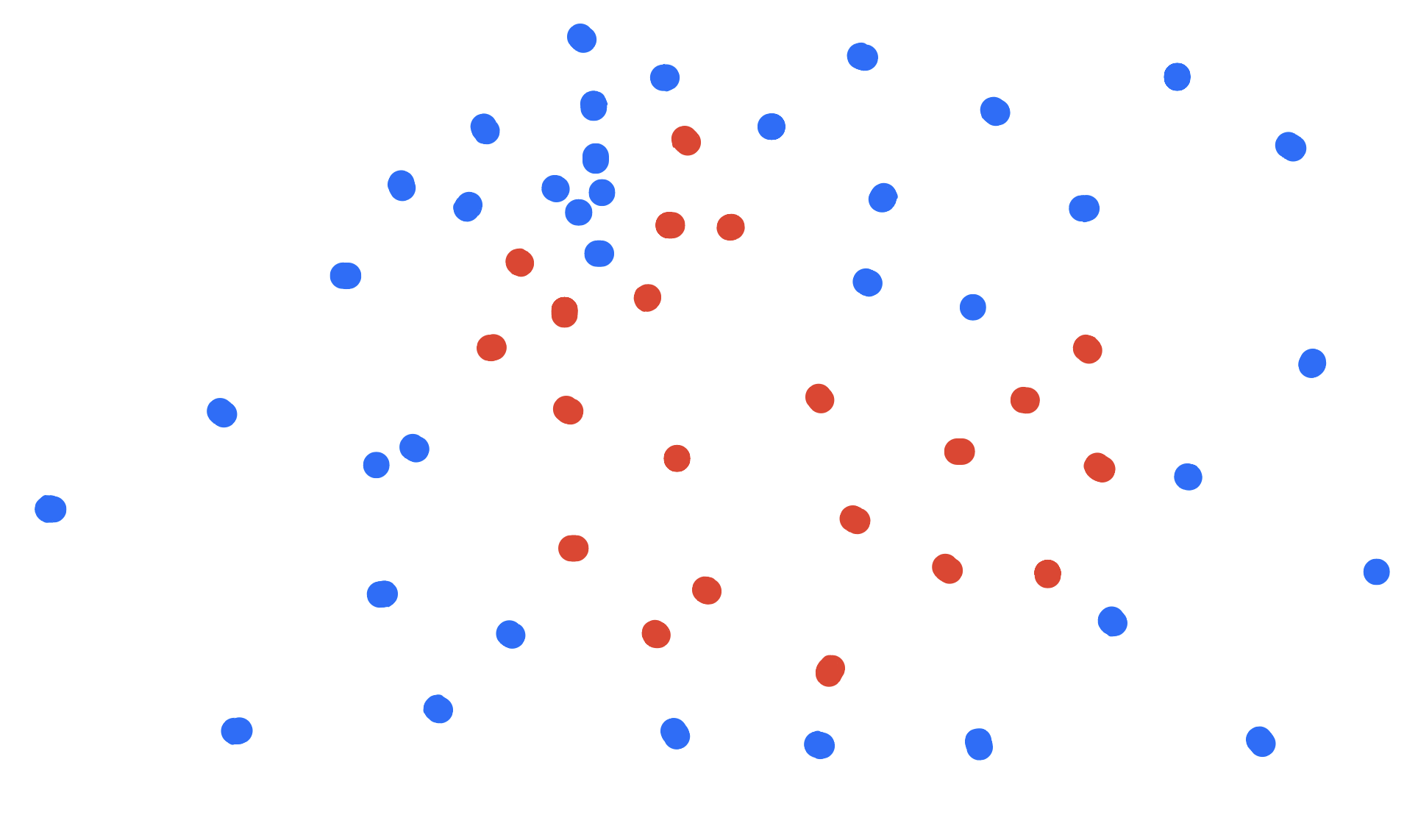}
  \]
Linear methods such as SVM can nevertheless still be used in many
cases, after application of a suitable {\it feature map}, namely a
(non-linear) transformation whose application on the data makes linear separation of classes possible:
\[
  \begin{array}{c}
    \includegraphics[width=5cm]{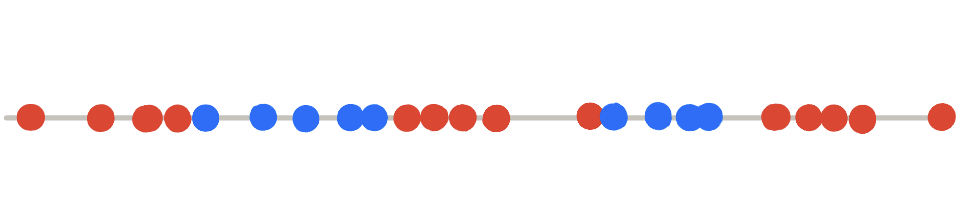}
  \end{array}
  \stackrel{\text{non-linear}}{\longto}
\begin{array}{c}  \includegraphics[width=4.5cm]{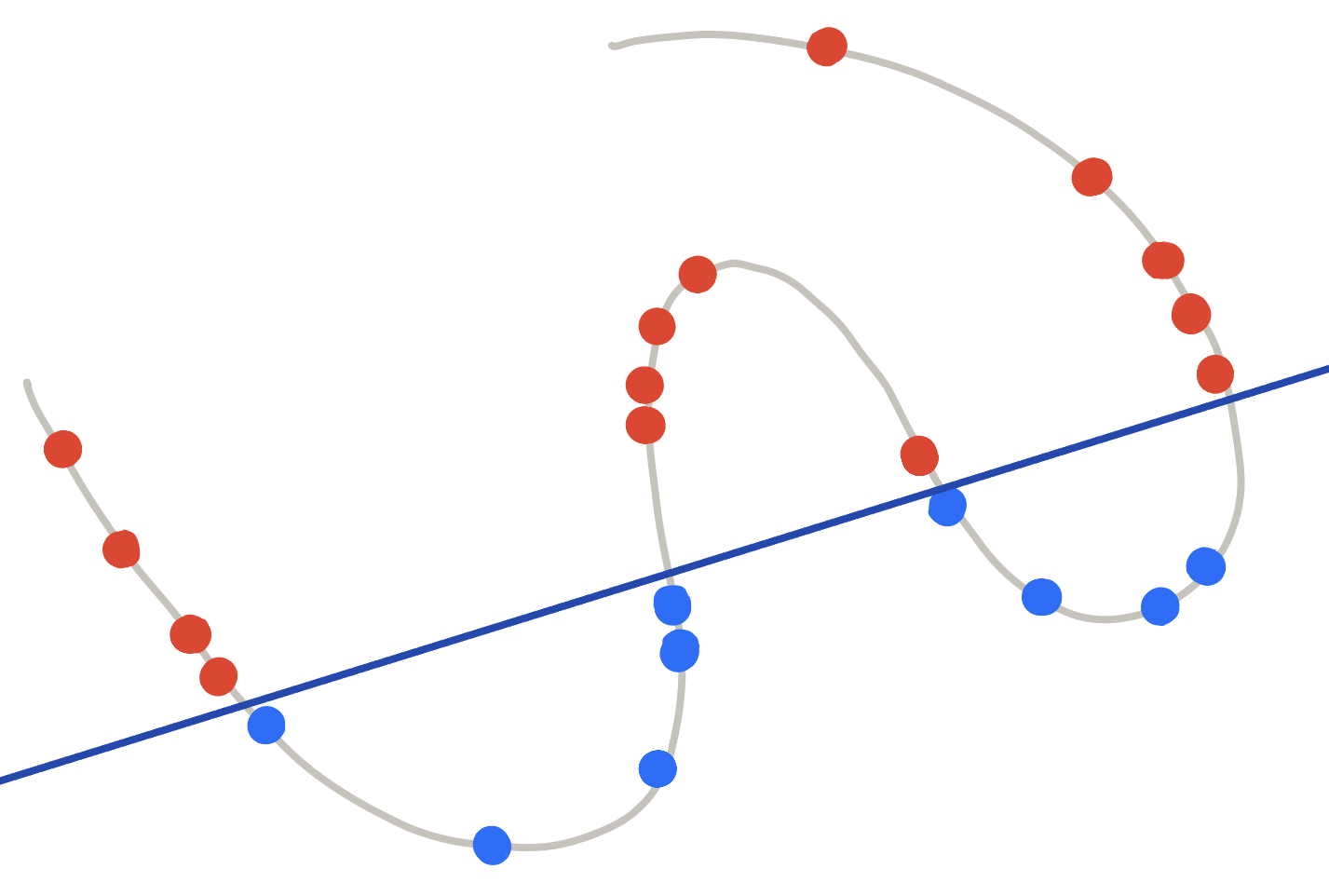} \end{array}
\]
It is on such more difficult learning tasks that deep learning can come into its own. The idea is that
  successive layers of the neural net transform the data gradually,
  eventually leading to an easier learning problem.\footnote{This idea seems to have been present in the
  machine learning literature for decades, see
  e.g. \cite{lecun1998gradient}. It is well explained in \cite[\S
  6]{goodfellow2016deep}.  For illustrations of this as well as the connection
to fundamental questions in topology, see the work of Olah
\cite{NNOlah}.}

In the standard setting of supervised learning, we assume the
existence of a function
\[
\phi : \RM^n \to \RM^m
\]
and know a (usually large) number of its values.
The task is to find a reasonable approximation of $\phi$,
given these known values. (The reader should keep in mind the
motivating problem of the previous
section, where one wants to learn a function
\[
  \phi : \RM^{28 \times 28} \to \RM^{10}
\]
giving the probabilities that a certain $28 \times 28$-pixel grayscale
image represents one of the 10 digits $0$, $1$, \dots, $9$.)

We fix a network architecture, which in our
simple setting of a vanilla neural net, means that we fix the
number of layers $\ell$ and  layer  dimensions $n_2, \dots, n_{\ell
  -1}$. We then build a neural net (see \S \ref{sec:nn}) which serves
as our function approximator:
\begin{equation} \label{eq:nn}
\phi_{\approx}:   \RM^{n} \stackrel{A_1}{\longto} \RM^{n_2}
  \stackrel{f}{\longto} \RM^{n_2} \stackrel{A_2}{\longto}
  \RM^n_3
\stackrel{f}{\longto} \RM^{n_3} \stackrel{A_3}{\longto}  \dots \stackrel{f}{\longto}  
\RM^{d_{n-1}} \stackrel{A_{\ell-1}}{\longto}\RM^{m} 
\end{equation}
To begin with, the affine linear maps $A_i$ are initialised via 
some (usually random) initialization scheme, and hence the function $\phi_{\approx}$ output by
our neural network will be random and have no relation to our target function
$\phi$. We then measure the distance between our function
$\phi_{\approx}$ and $\phi$ via some \emph{loss function} $L$. (For
example, $L$ might be the mean squared distance between the
values of $\phi$ and $\phi_{\approx}$.\footnote{There are many subtleties here, and a good choice of
  loss function is one of them. In my limited experience, neural
  networks do a lot better learning probability distributions than
  general functions. When learning probability
  distributions, \emph{cross entropy} \cite[\S 3.13]{goodfellow2016deep}
  is the loss function of
  choice.}) A crucial assumption is that this loss function is
differentiable
  in terms of the weights of our neural network. Finally, we perform
  gradient descent with respect to the loss function, in order to update
  the parameters in \eqref{eq:nn} to (hopefully) better and better
  approximate $\phi$.

  In order to get an intuitive picture of what is happening during
  training, let us assume that $m = 1$ (so we are trying to learn a
  scalar function), and that our activation functions are ReLU. Thus
  $\phi_{\approx}$ is the composition of affine linear and piecewise
  linear functions, and hence is piecewise linear. As with any
  piecewise linear function, we obtain a decomposition of $\RM^n$ 
  into polytopal regions
  \begin{equation} \label{eq:regions}
    \begin{array}{c}
\includegraphics[width=5cm]{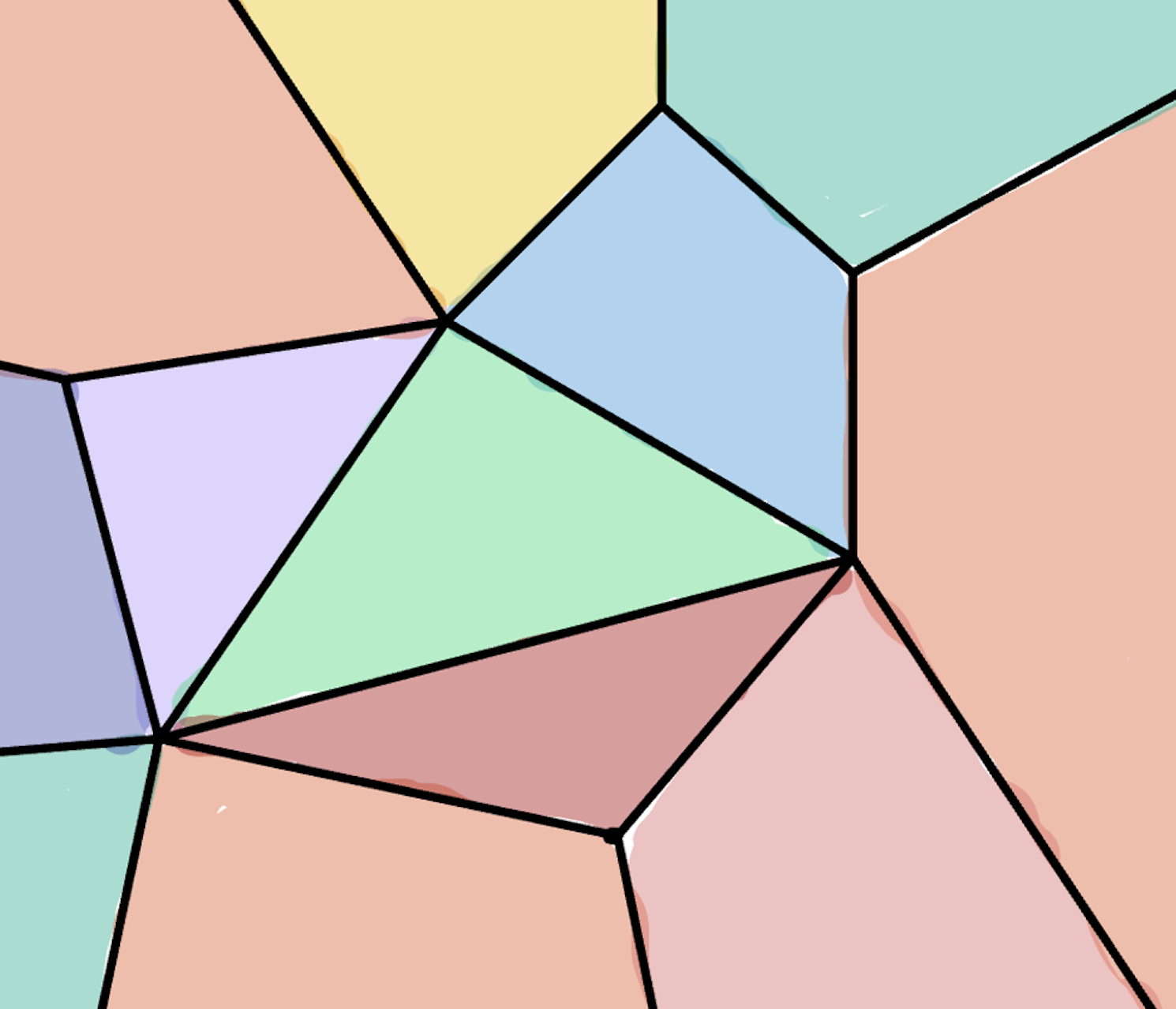} \end{array}
\end{equation}
such that $\phi_{\approx}$ is affine linear on each region. As training
progresses, the affine linear functions move in a way similar to
the learning of a line of best fit, but more complex since
the regions we are dealing with may also move, disappear or spawn new
regions.
  \begin{remark}
    For an excellent interactive animation of a simple neural network
    learning a classification task, the reader is urged to experiment with
    the Tensor Flow Playground \cite{tfp}. Karpathy's convolutional
    neural network demo \cite{karpathy} is also illustrative to
    play with.
  \end{remark}

  \begin{remark} Some remarks:
    \begin{enumerate}
    \item Typically, one splits the known values of $\phi$ into two
      disjoint sets, consisting of \emph{training data} and \emph{validation
        data}. Steps of gradient descent are only performed using
      the training data and the validation data allows us to periodically check whether
      our model is also making reasonable predictions at points not
      present in the training data (``validation error''). It is
      sometimes useful to have an additional set of \emph{test data}, completely unseen
      during training, which one can use to compute the performance of
      the trained model (``test error'').
    \item In most applications of machine learning, the training data
      is enormous and feeding it all through the neural network
      \eqref{eq:nn} in order to compute the loss function is unduly expensive. Thus one usually employs
      \emph{stochastic gradient descent}: at every step the gradient
      of the loss function is computed using a small random subset
      (a ``minibatch'') of the training data.
      \item Using a model with a small number of parameters (as
        traditionally done in statistics, and some ML methods) has
        advantages for interpretability and computation. It can also
        help avoid \emph{overfitting}, where the chosen
        predictor may fit the data set so closely it ends up fitting
        noise and fails to adequately capture the underlying data
        generating process. A simple example of a model with few
        parameters in a line of best fit. Deep learning is
      different, in that often there are enough parameters 
      to allow overfitting. What is surprising is that often neural
      nets generalize well
      (i.e. don't overfit) even though they could in
      principle (this is an enormous subject, see
      e.g. \cite{belkin2019reconciling}).
    \end{enumerate}
  \end{remark}
  
\section{Simple examples from pure mathematics}

It is important to keep in mind that the main
motivating applications for deep learning research are very different
from those arising in pure mathematics. For example, the
``recognize a hand-written digit'' function considered in the previous
two sections is rather different to the Riemann zeta
function\footnote{One should keep in mind that neural networks are
  universal approximators: a large enough neural network can
  approximate any continuous function accurately 
  \cite{wiki-universal}. However, in practice some functions are much
  more easily learnt than others.}!

This means that the mathematician wanting to use machine learning
should keep in mind that they are using tools designed 
for a very different purpose. The hype that ``neural nets can learn
anything'' also doesn't help. The following rules of thumb are useful
to keep in mind when selecting a problem for deep learning:
\begin{enumerate}
\item \emph{Noise stable.} Functions involved in image and speech
  recognition motivated much research in machine learning. These
  functions typically have very high-dimensional input (e.g. $\RM^{
    100
    \times 100}$ for a square $100 \times 100$ grayscale image) and are
  noise stable. For example, we can usually recognise an image or
  understand speech after the introduction of a lot of noise. Neural
  nets typically do poorly on functions which are very
  noise-sensitive.\footnote{This point should be read with some
    caution. For example, evaluation of board positions in Go is not a
    particularly noise-stable problem.}
\item \emph{High dimensional.} If one thinks of a neural network as a
  function approximator, it is a function approximator that comes
  into its own on high-dimensional input. These are the settings in
  which traditional techniques like Fourier series break
  down, due to the curse of dimensionality. Deep learning should be
  considered when the difficulty comes from the dimensionality, rather
  than from the inherent complexity of the function.
\item \emph{Unit cube.} Returning to our (unreliable) analogy with
  biological neural nets, one expects all charges occurring in the
  brain to belong to some fixed small interval. The same is true of
  artificial neural networks: they perform best when all real numbers
  encountered throughout the network from input to output belong to some
  bounded interval. Deep learning 
  packages are often written assuming that the inputs
  belong to the unit cube $[0,1]^{n} \subset \RM^n$.
\item \emph{Details matter.}
  Design choices like network architecture and size, initialization scheme, choice of
  learning rate (i.e. step size of gradient descent), choice of optimizer etc. matter
  enormously. It is also important how the inputs to the neural
  network are encoded as vectors in $\RM^n$ (the
  \emph{representation}).\footnote{It seems silly to have to write
    that details matter in any technical subject. However many
    people I have spoken to are under the false impression that one model
    works for everything, and that training happens ``out of the
    box'' and is easy. For an excellent and honest summary by an expert of the difficulties 
    encountered when training large models, see \cite{karpathy2}.}
  Overcoming these difficulties is best done with a collaborator
  who has experience in deep learning research and implementation.
\end{enumerate}
With these rules of thumb in mind we will now discuss three examples
in pure mathematics.

\subsection{Learning the parity bit} Consider the parity bit function
\begin{align*}
  \sigma : \{ 0, 1 \}^m & \to \{0,1\}  \\
  (x_i) & \mapsto \sum_{i=1}^m x_i \mod 2.
\end{align*}
We might be tempted to use a neural network to try to learn a function 
\[
\sigma_{\approx} : \RM^m \to \RM
\]
which agrees with $\sigma$ under the natural embedding $\{ 0, 1 \}^m
\subset \RM^m$.

This is a classic problem in machine learning
\cite[I~\S~3.1]{minsky2017perceptrons}. It generalizes the problem of learning
the XOR function (the case $m = 2$), which is one of the simplest problems which cannot
be learned without non-linearities. There exist elegant neural
networks extending $\sigma$ to the unit cube, and given a large proportion (e.g. 50\%) of the
set $\{ 0, 1 \}^m$ a neural network can be trained to express $\sigma$
\cite[pp. 14-16]{rumelhart1985learning}.  However, given only a small
proportion of the values of $\sigma$ (e.g. $10\%$ for $m = 10$) a
vanilla neural network will not reliably generalize to all values of
$\sigma$ (for experiments, see \cite[`Playing with parity']{MLWM}).

The issue here is that $\sigma$ is highly noise sensitive. (Indeed, 
$\sigma$ is precisely the checksum of signal processing!) This is an
important example to keep in mind, as many simple functions in pure
mathematics resemble $\sigma$. For example, see \cite[Week 2]{MLWM}
where we attempt (without much luck!) to train a neural network to learn the M\"obius
function from number theory.

\subsection{Learning descent sets} Consider the symmetric group
$\Sigma_n$
consisting of all permutations of ${1, 2, \dots, n}$. Given a
permutation we can consider its \emph{left} and \emph{right descent sets}:
\begin{align}
  \LC(x) &= \{ 1 \le i < n \;  | \; x^{-1}(i) > x^{-1}(i+1) \}, \\
\RC(x) &= \{ 1 \le i < n \; | \; x(i) > x(i+1) \}. \label{eq:R}
\end{align}
Obviously, $\LC(x^{-1}) = \RC(x)$ and $\RC(x^{-1}) = \LC(x)$. The left
an right descent sets are important invariants of a permutation.

It is interesting to see whether a neural network can be trained to
learn the left and right descent sets. In other words, we would like
to train a neural network
\[
\phi_{\approx} : \RM^n \to \RM^{n-1}
\]
which given the vector $(x(1), x(2), \dots, x(n))$  returns a sequence
of $n-1$ probabilities giving whether or not $1 \le i < n$ belongs to the left
(resp. right) descent set.

This example is interesting in that \eqref{eq:R} implies that the right descent set can be
predicted perfectly with a single linear layer. More precisely, if we consider
\begin{align*}
  \gamma : \RM^n &\to \RM^{n-1} \\
  (v_1, \dots, v_n) &\mapsto (v_1 - v_2, v_2 - v_3, \dots, v_{n-1} - v_n)
\end{align*}
then the $i^{th}$ coordinate of $\gamma$ evaluated on a permutation $(x(1), \dots, x(n))$
is positive if and only if $i \in
\RC(x)$. On the other hand, it seems much harder to handcraft a
neural network which extracts the left descent set from $(x(1), \dots, x(n))$.

This might lead us to guess that a neural network will have a much easier time
learning the right descent set than the left descent set. This turns
out to be the case, and the difference is dramatic: a vanilla neural
network with two hidden layers of dimensions 500 and 100 learns to predict
right descent sets for $n = 35$ with high accuracy after a few
seconds. Whereas, the same network struggles to get even a single
correct answer for the left descent set, after significant training!\footnote{Decreasing $n$ and
  allowing longer training suggests that the network can learn the
  left descent set, however it is \emph{much} harder.} It is 
striking that using permutation matrices as inputs rather than the
vectors $(x(1), \dots, x(n))$ gives perfect symmetry in training
between left and right.\footnote{For a colab containing all of these
experiments, see \cite[Classifying descent sets in $S_n$]{MLWM}.}
The issue here is the \emph{representation}: how the model receives
its input can have a dramatic effect on model performance.
  
\subsection{Transformers and linear algebra} Our final example is much
more sophisticated, and illustrates how important the choice of training data can
be. It also shows how surprising the results of training large
neural networks can be.

A \emph{transformer} is a neural network architecture which first emerged in
machine translation \cite{vaswani2017attention}. We will not go into
any detail about the transformer architecture here, except to say that
it is well-suited to tasks where the input and output are sequences of
tokens (``sequence to sequence'' tasks):
\[
  \begin{tikzpicture}[scale=.8]
     \draw (-3.5,-.5) rectangle (-2.5,.5);
     \draw (-2.5,-.5) rectangle (-1.5,.5);
     \draw (-1.5,-.5) rectangle (-0.5,.5);
     \node (x) at (-3,0) {$\begin{array}{c} x \end{array}$};
     \node (y) at (-2,0) {$\begin{array}{c} y \end{array}$};
     \node (z) at (-1,0) {$\begin{array}{c} z \end{array}$};
     \draw[->] (0,0) to (1,0);
     \draw (1.5,-.5) rectangle (4.5,.5);
     \node (text) at (3,0) {transformer};
     \draw[->] (5,0) to (6,0);
     \draw (6.5,-.5) rectangle (7.5,.5);
     \draw (7.5,-.5) rectangle (8.5,.5);
     \draw (8.5,-.5) rectangle (9.5,.5);
          \node (a) at (7,0) {$\begin{array}{c} a \end{array}$};
     \node (b) at (8,0) {$\begin{array}{c} b \end{array}$};
     \node (c) at (9,0) {$\begin{array}{c} c \end{array}$};
   \end{tikzpicture}
  \]
More precisely, the input sequence (``xyz'') determines a probability
distribution over all tokens. We then sample from this distribution to
obtain the first token (``a''). Now the input and sequence sampled so far
(``xyz'' + ``a'') provides a new distribution over tokens, from which
we sample our second token (``b''), etc.
  
In a recent work \cite{Charton}  Charton trains a transformer to
perform various tasks in linear
algebra: matrix transposition, matrix addition, matrix multiplication, determination of
eigenvalues, determination of eigenvectors etc. For example, the
eigenvalue task is regarded as the 
``translation'':
\[
    \begin{tikzpicture}[scale=.75]
  \node (l) at (-3,0) {$\begin{array}{c} \text{real $5 \times 5$-symmetric matrix} \\ M =
                   (m_{11}, m_{12}, m_{13}, \dots, m_{55}) \end{array}$};
     \draw[->] (l) to (1,0);
     \draw (1.5,-.5) rectangle (4.5,.5);
     \node (text) at (3,0) {transformer};
     \node (r) at (8,0) {$\begin{array}{c} \text{list of eigenvalues} \\
                       \lambda_1 \ge \lambda_2 \ge \dots \ge \lambda_5.
                         \end{array}$};
                            \draw[->] (5,0) to (r);
         \end{tikzpicture}
       \]
Charton considers real symmetric matrices, all of whose entries are signed
floating point numbers with three significant figures and exponent
lying between $-100$ and $100$.\footnote{Charton considers various
  encodings of these numbers via sequences of tokens of various
  lengths, see \cite{Charton}.} The
transformer obtains impressive accuracy on most linear algebra tasks. 
What is remarkable is that for the transformer 
the entries of the matrix (e.g. 3.14, -27.8, 0.000132, \dots) are simply tokens---the
transformer doesn't ``know'' that 3.14 is close to 3.13, or that both
are positive; it doesn't
even ``know'' that its tokens represent numbers!

Another remarkable aspect of this work concerns generalization. A
model trained on Wigner matrices (e.g. entries sampled
uniformly from $[-10,10]$) does not generalize well at all to matrices
with positive eigenvalues. On the other hand, 
a model trained on matrices with eigenvalues sampled from a Laplace
distribution (which has heavy
tails) does generalize to matrices whose eigenvalues are all positive,
even though it has not seen a single such matrix during training! The
interested reader is referred to Charton's paper 
\cite{Charton} (in particular Table 12) and his lecture on youtube
\cite{ChartonYoutube}.
       
\section{Examples from research mathematics}

We now turn to some examples where deep learning has been used in pure
mathematics research.

\subsection{Counter-examples in combinatorics} One can dream that
deep learning might one day provide a mathematician's ``bicycle for
the mind'': an easy to use and flexible framework for
exploring possibilities and potential counter-examples. (I have
certainly lost many days trying to prove a statement that
turned out to be false, with the counter-example lying just beyond my
mental horizon.)

We are certainly not there yet, but the closest we have come to
witnessing such a framework is provided in the work of Adam Wagner
\cite{wagner2021constructions}. He focuses on conjectures of the
form: over all combinatorial structures $X$, an associated numerical
quantity $Z$ is bounded by $B$. He considers situations where there is some simple
recipe for generating objects in $X$, and that the numerical 
quantity $Z$ is efficiently computable.

For example, a
conjecture in graph theory states that for any connected graph $G$ on
$n \ge 3$ vertices, with largest eigenvalue $\lambda$ and matching number $\mu$
we have
\begin{equation} \label{eq:reward}
\lambda + \mu - \sqrt{n-1} - 1 \ge 0.
\end{equation}
(It is not important for this discussion to know what the matching
number or largest eigenvalue are!)

  \begin{figure}[t] \label{fig:ev}
    \caption{The evolution of graphs towards Wagner's
      counter-example, from \cite{wagner2021constructions}.}
\[    \begin{array}{c}
      \\
      \includegraphics[width=12cm]{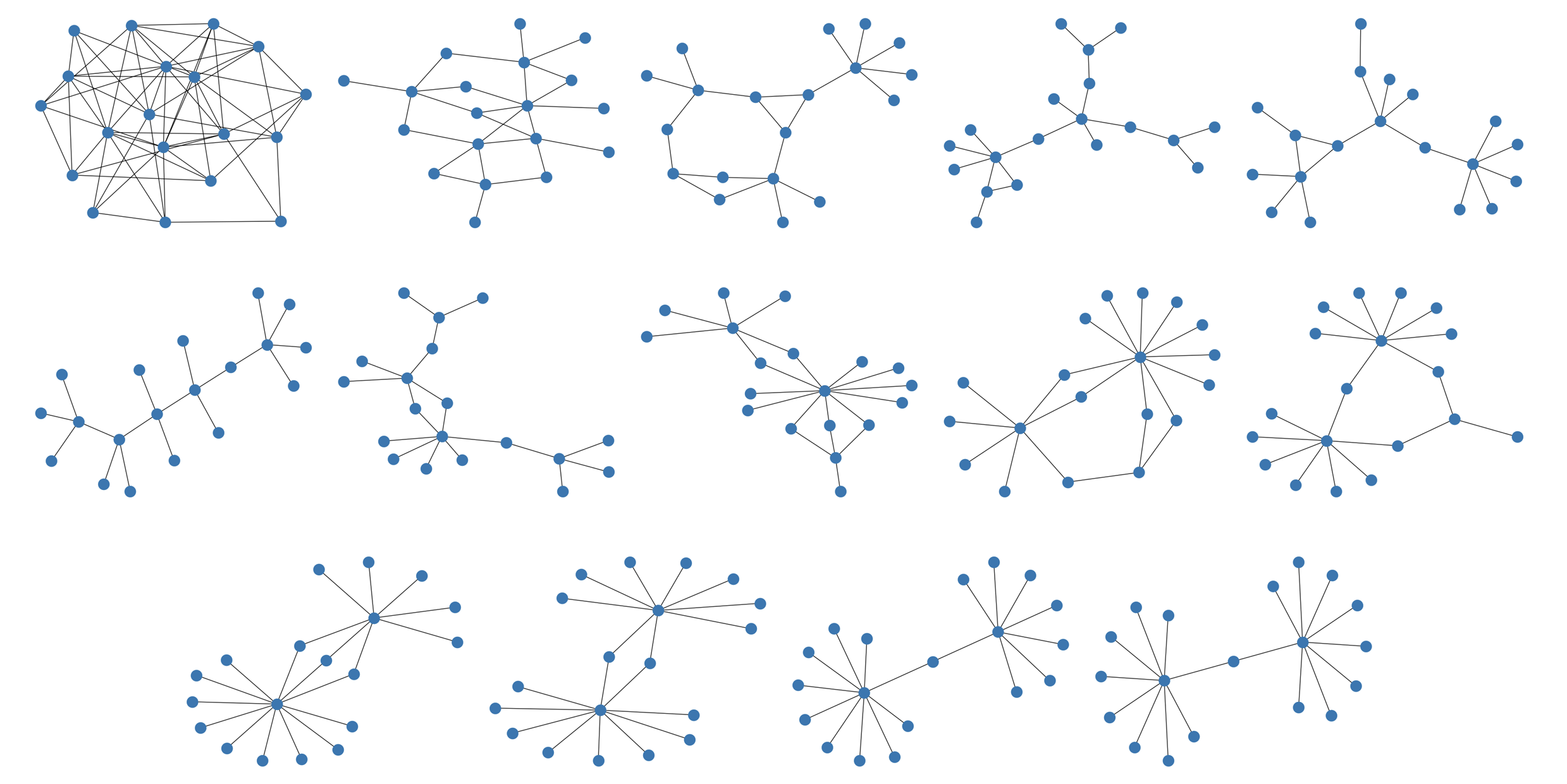}
      \end{array} \]
    \end{figure}

Wagner fixes an enumeration $e_1, e_2, \dots $ of the edges $E$ in a
complete graph on $n$-vertices. Graphs are generated by playing a
single player game: the player is offered $e_1$, $e_2$ etc. and
decides at each point whether to accept or reject the edge, the goal
being to minimize \eqref{eq:reward}. A move in the game is 
given by a $01$-vector indicating edges that have been taken so far,
together with a vector indicating which edge is under
consideration. For example, when $n = 4$ the pair $((1,0,1,1,0,0), (0,0,0,0,1,0))$
indicates that edge number 5 is under consideration, and that edges
$1$, $3$ and $4$ have already been selected, and 2 rejected. Moves are
sampled according to a neural network
\begin{equation} \label{eq:Wnn}
\mu :\RM^E \oplus \RM^E \to \RM,
\end{equation}
which (after application of sigmoid) gives the probability
that we should take the edge under consideration.

Wagner then employs the \emph{cross entropy method} to gradually
train the neural network. A fixed (and large) number of graphs are
sampled according to the neural network \eqref{eq:Wnn}. Then a fixed
percentage (say 10\%) of the games resulting in the smallest values of
the LHS of \eqref{eq:reward} are used as training data to update the neural
network \eqref{eq:Wnn}. (That is, we tweak the
weights of the neural network to make decisions that result in graphs
that are as close as possible to providing a counter-example to
\eqref{eq:Wnn}.) We then repeat. This
method eventually finds a counter-example to \eqref{eq:reward} on 19
vertices. The evolution of graphs sampled from the neural network is
shown in Figure \ref{fig:ev}---note how the neural network learns
quickly that tree-like graphs do best. Exactly the same method works to discover
counter-examples to several other conjectures in combinatorics, see
\cite{wagner2021constructions}.

\subsection{Conjecture generation}

    \begin{figure}[b]
\[
\begin{array}{c}
    \includegraphics[scale=0.25]{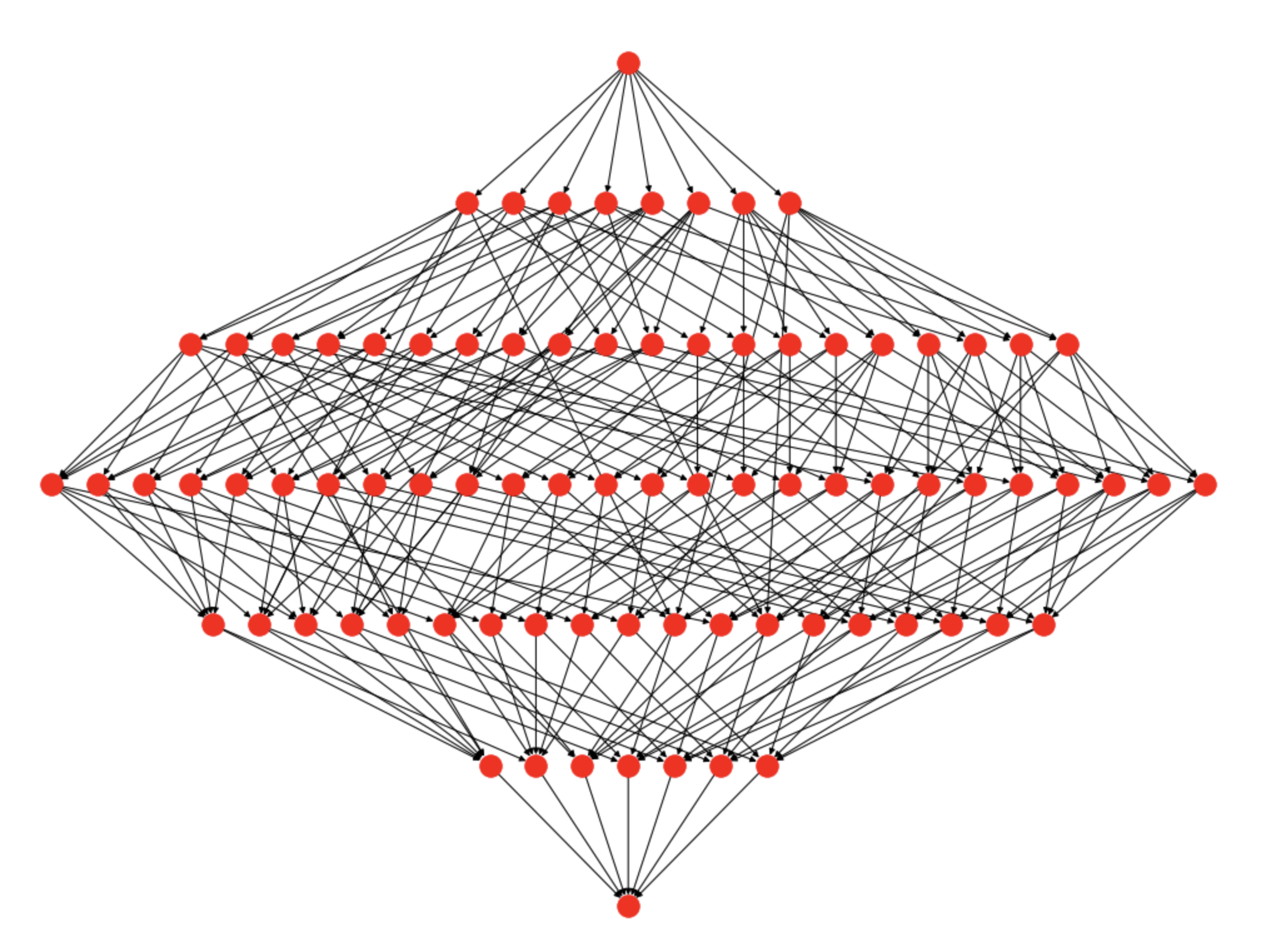}
\end{array}
\leftrightarrow 1 + 3q + q^2
\]
  \caption{Bruhat interval and Kazhdan-Lusztig polynomial for the pair
    of permutations $x =
    (1,3,2,5,4,6)$ and $y=(3,4,5,6,1,2)$ in $\Sigma_6$, from \cite{ci}}.
  \label{fig:KL131}
\end{figure}

    The \emph{combinatorial invariance conjecture} is a conjecture in
    representation theory which was proposed by Lusztig and Dyer in
    the early 1980s \cite{brenti}. To any pair of permutations $x, y
    \in \Sigma_n$ in the symmetric group
    one may associate two objects: the \emph{Bruhat graph} (a directed
    graph); and the \emph{Kazhdan-Lusztig polynomial} (a polynomial in
    $q$), see Figure \ref{fig:KL131} for an example of both. The
    conjecture states that an isomorphism between Bruhat 
    graphs implies equality between Kazhdan-Lusztig polynomials. A
    more optimistic version of this conjecture asks for a recipe which
    computes the Kazhdan-Lusztig polynomial from the Bruhat graph. One
    interesting aspect of this conjecture is that it is (to the best of my
    knowledge) a conjecture born of pure empiricism.

    For the Bruhat graph, the definition is simple, but the resulting
    graph is complicated. On the other
    hand, the definition of the Kazhdan-Lusztig polynomial is
    complicated, however the resulting polynomial is simple. Thus,
    there is at least a passing resemblance to traditional
    applications of machine learning, where a simple judgement (e.g. ``it's
    a cat'') is made from complicated input (e.g. an array of
    pixels).

It is natural to use neural networks as a testing ground for
this conjecture: if a neural network can easily predict the
Kazhdan-Lusztig polynomial from the Bruhat graph, perhaps we can too!
We trained a neural network to predict Kazhdan-Lusztig
polynomials from the Bruhat graph. We used a neural network
architecture known as a graph neural network, and trained the neural
network to predict a probability distribution on the coefficients of
$q$, $q^2$, $q^3$ and $q^4$.\footnote{The coefficient of $q^0$ is known to
always equal $1$. In our training sets no coefficients of $q^5$ or
higher occur.} The neural network was trained on $\approx 20\;000$ Bruhat
graphs, and achieved very high accuracy ($\approx 98\%$) after less
than a day's training. This provides reasonable evidence that there is
\emph{some} way of reliably guessing the Kazhdan-Lusztig polynomial
from the Bruhat graph.

\begin{figure}
\[
\begin{array}{c}
    \includegraphics[scale=0.15]{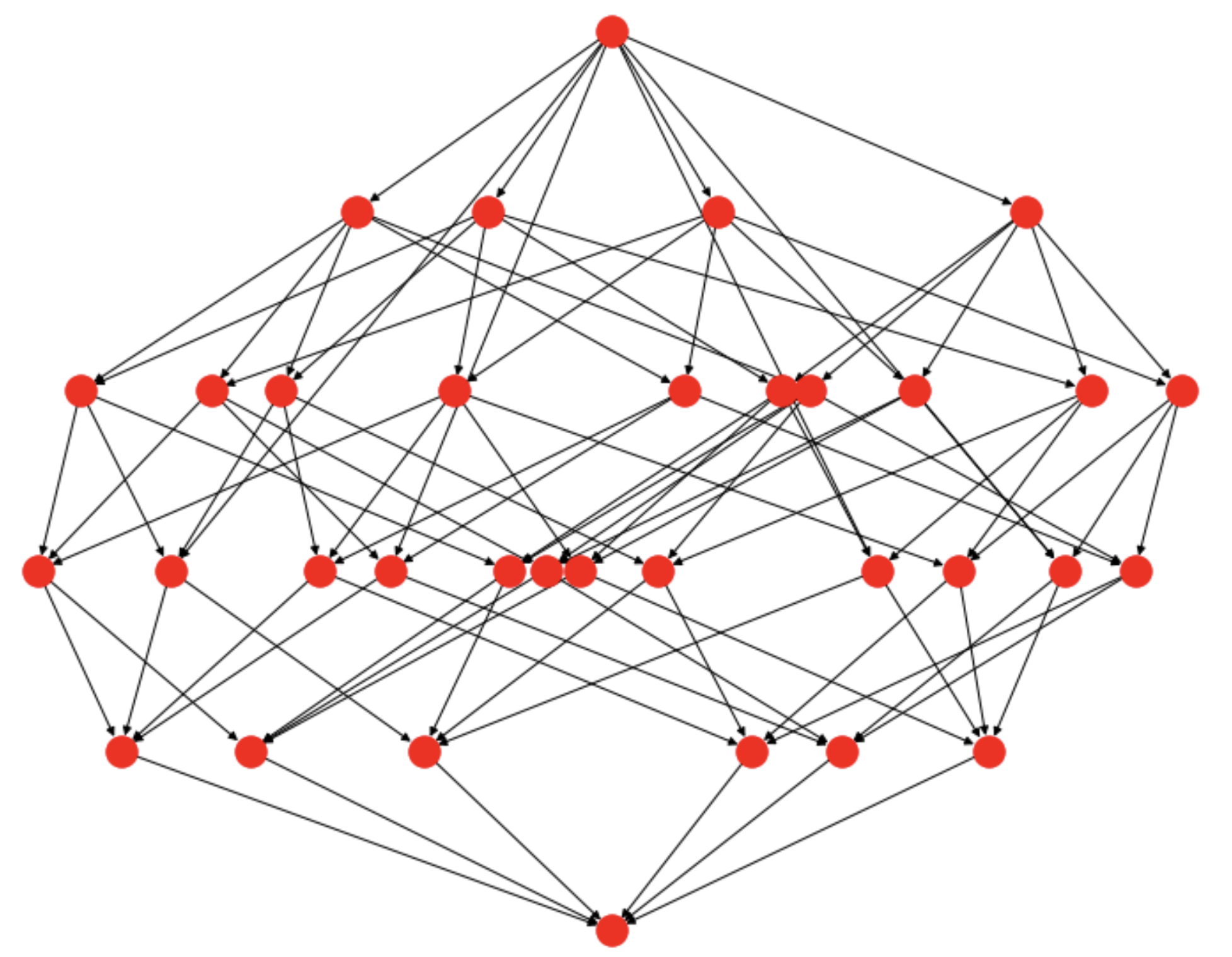}
\end{array}
\begin{array}{c}
    \includegraphics[scale=0.15]{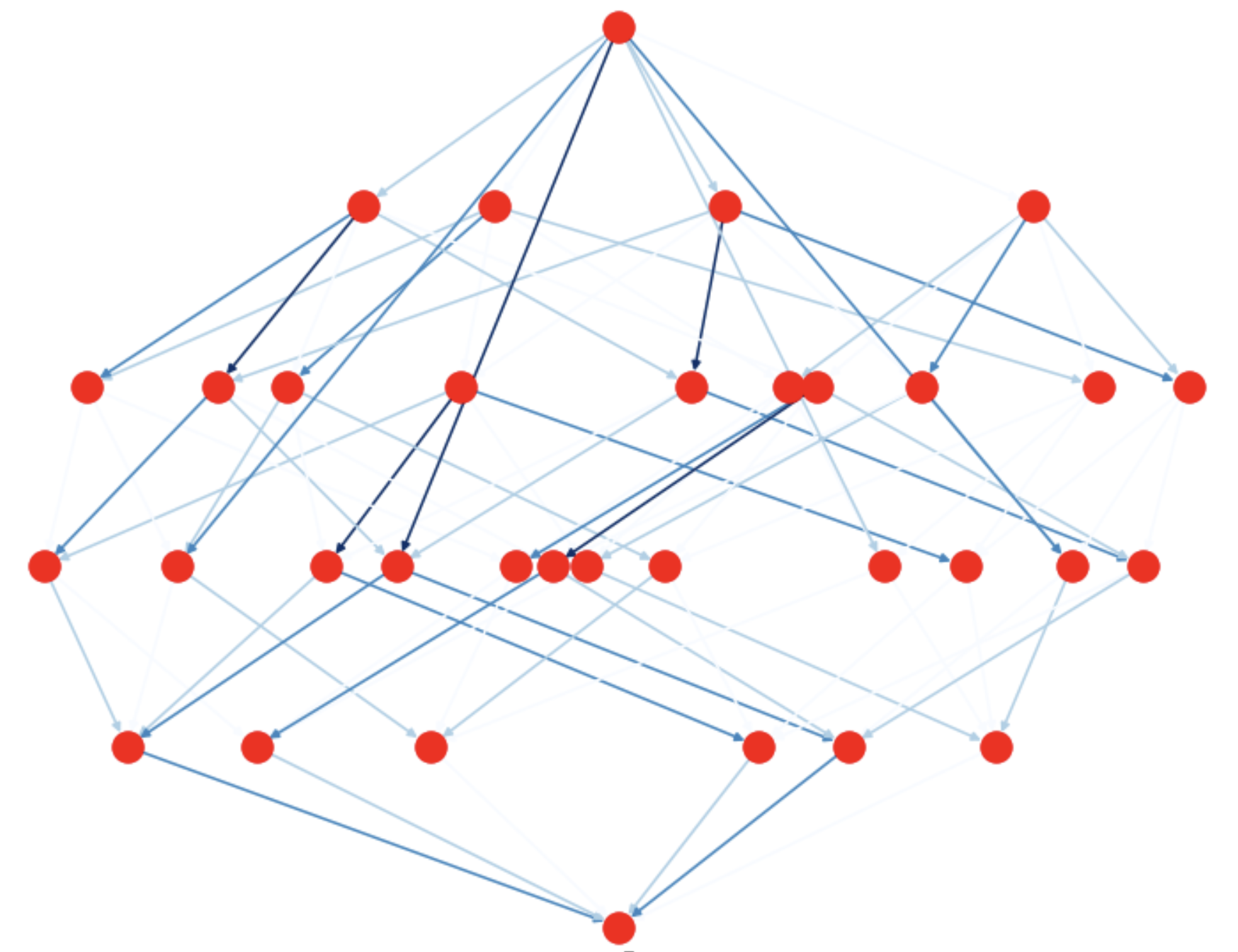}
\end{array}\]
\caption{Bruhat interval pre and post saliency analysis.}
  \label{fig:saliency}
\end{figure}

It is notoriously difficult to go from a trained neural network to
some kind of human understanding. One technique to do so is known as
\emph{saliency analysis}. Recall that neural
networks often learn a piecewise linear function, and hence one can
take derivatives of the learned function to try to learn which inputs
have the most influence on a given output.\footnote{This technique is 
often called ``vanilla gradient'' in the literature. Apparently it is very
brittle in real-world applications.} In our example, saliency analysis
provided subgraphs of the original Bruhat graph which appeared to have
remarkable ``hypercube'' like structure (see Figure
\ref{fig:saliency} and \cite[Figure 5a]{davies2021advancing}). After
considerable work this eventually led to a
conjecture \cite{ci} which would settle the combinatorial invariance
conjecture for symmetric groups if proven, and has stimulated research on this
problem from pure mathematicians \cite{gurevich2023parabolic,Gaetz,Gaetz2,brentim}.

In a parallel development, Davies, Juh{\'a}sz, Lackenby and Tomasev
were able to use saliency analysis to discover a new relationship
between the signature and hyperbolic invariants of knots
\cite{DJLT}. The machine learning background of both works is
explained in \cite{davies2021advancing}. It would be very interesting
to find further examples where saliency leads to new conjectures and
theorems.

\subsection{Guiding calculation} Another area where deep learning has
promise to impact mathematics is in the guiding of calculation. In many settings a computation can be done in
many ways. Any choice will lead to a correct outcome, but
choices may
drastically effect the length of the computation. It is interesting to
apply deep learning in these settings, as false steps (which
deep learning models are bound to make) effects efficiency but not
accuracy.

Over the last three years there have been several examples of such
applications. In \cite{peifer2020learning}, the authors use a machine
learning algorithm to guide selection strategies in Buchberger's
algorithm, which is a central algorithm in the theory of Gr\"obner
bases in polynomial rings. In \cite{simpson2021learning}, Simpson
uses deep neural networks to simplify proofs in the classification of
nilpotent semi-groups. In \cite{heal2022deep}, the authors use a deep
neural network to predict computation times of period matrices, and
use it to more efficiently compute the periods of certain hypersurfaces in
projective space.

\subsection{Prediction} Due to limitations of space, we cannot begin to
survey all the work done in this infant subject. In particular, there
has been much work (see
e.g. \cite{bao2021polytopes,brodie2020machine}) training 
neural networks to predict difficult quantities in mathematics
(e.g. volumes of polytopes, line bundle cohomology,\dots).

\section{Conclusion}

The use of deep learning in pure mathematics is in its infancy. The
tools of machine learning are flexible and powerful, but need
expertise and experience to use. One should not expect things to work
``out of the box''. Deep learning has found applications
in several branches of pure mathematics including combinatorics,
representation theory, topology and algebraic geometry. Applications
so far support the thesis that deep learning most usefully aids
the more intuitive (``system 1'') parts of the mathematical process:
spotting patterns, deciding where counter-examples might lie, choosing
which part of a calculation to do next. However, the possibilities do seem
endless, and only time will tell.

\bibliographystyle{myalpha}
\bibliography{gen}

  \end{document}